\documentclass[a4paper,11pt]{amsart}

\usepackage[utf8]{inputenc}
\usepackage[T1]{fontenc}
\usepackage{lmodern}
\usepackage[a4paper,margin=1in]{geometry}
\usepackage{microtype}
\usepackage{mathtools,amssymb,amsfonts}
\usepackage{enumitem}
\usepackage{hyperref}
\hypersetup{colorlinks=true,citecolor=blue,linkcolor=blue,urlcolor=blue}

\numberwithin{equation}{section}
\newtheorem{theorem}{Theorem}[section]
\newtheorem{corollary}[theorem]{Corollary}
\newtheorem{lemma}[theorem]{Lemma}
\newtheorem{proposition}[theorem]{Proposition}
\newtheorem{definition}[theorem]{Definition}
\newtheorem{remark}[theorem]{Remark}
\newtheorem{example}[theorem]{Example}

\newcommand{\OK}{\mathcal O_K}
\newcommand{\Z}{\mathbb Z}
\newcommand{\Q}{\mathbb Q}
\newcommand{\F}{\mathbb F}
\newcommand{\A}{\mathcal A}
\newcommand{\E}{\mathcal E}
\newcommand{\M}{\mathcal M}
\newcommand{\Ssq}{\mathcal S}
\newcommand{\vp}{\nu_p}
\DeclareMathOperator{\rad}{rad}
\DeclareMathOperator{\disc}{disc}

\title[$\alpha$-monogeneity of pure number fields]{$\alpha$-monogeneity of pure number fields:\\ criterion and density}

\author{Khai-Hoan Nguyen-Dang}
\address{Morningside Center of Mathematics, Chinese Academy of Sciences, No.\ 55, Zhongguancun East Road, Beijing 100190, China}
\email{khaihoann@gmail.com}

\author{Thai-Hung Nguyen}
\address{Department of Mathematics, Ho Chi Minh City University of Education, 280 An Duong Vuong Street, Cho Quan Ward, Ho Chi Minh City, Vietnam}
\email{thaihungspt2003@gmail.com}


\date{\today}

\begin{document}

\begin{abstract}
Let \(n\ge 2\), let \(m\in\Z\setminus\{0\}\), and let \(K=\Q(\alpha)\), where \(\alpha^n=m\) and \(X^n-m\) is irreducible over \(\Q\). We study when the natural order \(\Z[\alpha]\) is the full ring of integers \(\OK\). For the pure family \(X^n-m\), we give a short proof, using only Dedekind's index criterion, of the equivalence
\[
\OK=\Z[\alpha]
\iff
\text{\(m\) is square-free and \(\nu_p(m^p-m)=1\) for every prime \(p\mid n\).}
\]
Equivalently, the prime support of \([\OK:\Z[\alpha]]\) is
\[
\{p:\nu_p(m)\ge 2\}\cup \{p\mid n:\nu_p(m^p-m)\ge 2\}.
\]
We then compute the natural density of the corresponding parameters in the one-parameter family \(X^n-m\):
\[
\delta_n=\frac{6}{\pi^2}\prod_{p\mid n}\frac{p}{p+1}.
\]
We also give an arithmetic-progression refinement, a density-theoretic independence statement for the local obstruction sets at primes dividing \(n\), and discriminant-ordered counts of the corresponding fields.
\end{abstract}

\maketitle

\section{Introduction}\label{sec:intro}

Let \(n\ge 2\) be fixed. For a nonzero integer \(m\) such that \(X^n-m\) is irreducible over \(\Q\), choose a root \(\alpha_m\) of \(X^n-m\), and put \(K_m=\Q(\alpha_m)\). The order \(\A_m=\Z[\alpha_m] \) is the natural order attached to the defining equation \(X^n-m\). We study the following question.
\begin{quote}
For which parameters \(m\) is \(\A_m\) the maximal order of \(K_m\), that is, when is \(\A_m=\mathcal O_{K_m}\)?
\end{quote}
When \(m>1\) is square-free, \(X^n-m\) is Eisenstein at every prime dividing \(m\), hence irreducible, and the parameter space is the pure family of degree \(n\) number fields.

This question is related to, but distinct from, field monogeneity. The field \(K_m\) is monogenic if \(\mathcal O_{K_m}=\Z[\theta]\) for some algebraic integer \(\theta\in K_m\). The equality \(\A_m=\mathcal O_{K_m}\) asks for a power basis generated by the distinguished root \(\alpha_m\). Thus \(\A_m=\mathcal O_{K_m}\) implies field monogeneity, but the converse may fail for this specified generator; already in degree two the two conditions differ on a positive-density set of square-free parameters.

Our first main theorem gives the density of the set of parameters for which the natural order is maximal. Let
\[
\M_n^+(T)=\#\{1\le m\le T: T^n-m\text{ is irreducible over }\Q\text{ and }\Z[\alpha_m]=\mathcal O_{K_m}\}.
\]

\begin{theorem}\label{thm:intro-density}
For every fixed \(n\ge 2\),
\[
\M_n^+(X)=\delta_n X+O_n(X^{1/2})\qquad (X\to\infty),
\qquad
\delta_n=\frac{6}{\pi^2}\prod_{p\mid n}\frac{p}{p+1}.
\]
Equivalently, relative to the square-free positive integers, the density of parameters \(m\) for which \(\A_m\) is the maximal order of \(K_m\) is \(\prod_{p\mid n}\frac{p}{p+1}.\) 
\end{theorem}

The density theorem rests on a purely local test for the index of \(\Z[\alpha]\), as follows. Throughout the paper, \(\nu_p\) denotes the \(p\)-adic valuation on \(\Z\), and we use the convention \(\nu_p(0)=+\infty\).

\begin{theorem}\label{thm:intro-criterion}
Let \(m\in\Z\setminus\{0\}\) be such that \(X^n-m\) is irreducible over \(\Q\), and let \(K=\Q(\alpha)\) with \(\alpha^n=m\). Then
\[
\Z[\alpha]=\mathcal O_K
\iff
\text{\(m\) is square-free and \(\vp(m^p-m)=1\) for every prime \(p\mid n\).}
\]
Here and below, an integer \(m\ne0\) is square-free if no prime square divides \(m\), equivalently \(\mu^2(|m|)=1\).
\end{theorem}

The statement of Theorem~\ref{thm:intro-criterion} is a specialization to the pure family of Smith's relative criterion for radical extensions \cite{SmithRadical2021}. Our contribution here is different: we record a short pure-case proof using only Dedekind's index criterion, and then use that criterion to obtain an exact density theory in the one-parameter family \(X^n-m\).

The preceding theorem concerns the distinguished generator \(\alpha_m\). To avoid any possible conflict with the usual notion of field monogeneity, we record the precise comparison needed later.

\begin{corollary}\label{cor:intro-comparison}
Let \(m>1\) run through the square-free integers.
\begin{enumerate}[label=\textup{(\arabic*)}]
\item If \(n=2\), then every \(K_m\) is monogenic, but \(\A_m\) is maximal if and only if \(m\not\equiv 1\pmod 4\).
\item If \(n=3\), then the set of square-free positive \(m\) for which \(K_m\) is monogenic but \(\A_m\ne\mathcal O_{K_m}\) has relative density zero inside the square-free positive integers. In particular, counted as a subset of the positive integers, the monogenic pure cubic parameters have density \(\delta_3=\frac{9}{2\pi^2}\).
\item If \(n\ge4\), then the set of square-free positive \(m\) for which \(K_m\) is monogenic but \(\A_m\ne\mathcal O_{K_m}\) has relative density zero inside the square-free positive integers. In particular, counted as a subset of the positive integers, the monogenic pure degree-\(n\) parameters have density \(\delta_n\).
\end{enumerate}
\end{corollary}

Part~(1) is classical, see for example \cite[Chapter 1]{NeukirchANT}. Part~(2) uses the work of Aygin--Nguyen on monogenic pure cubics \cite{AyginNguyen2021,AyginNguyenCorrigendum2023}, while part~(3) follows from the first author's work \cite[Theorem~5.4]{ND26}. Thus, except in degree two, Theorem~\ref{thm:intro-density} gives the absolute positive-parameter density of genuine field monogeneity in the pure family.

For parameters \(m\) satisfying \(\Z[\alpha_m]=\mathcal O_{K_m}\), the index
\([\mathcal O_{K_m}:\Z[\alpha_m]]\) is equal to \(1\). Hence the field
discriminant is the polynomial discriminant. This allows the parameter count to be translated directly into a discriminant-ordered field count.

\begin{theorem}\label{thm:intro-disc}
Let \(N_n(Y)\) denote the number of isomorphism classes of fields \(K=\Q(\alpha)\) for which \(\alpha^n=m\) for some \(m\in\Z\setminus\{0\}\), the polynomial \(X^n-m\) is irreducible over \(\Q\), \(\Z[\alpha]=\mathcal O_K\), and \(|\disc(K)|\le Y\). Put
\[
X=\left(\frac{Y}{n^n}\right)^{1/(n-1)}.
\]
Then, as \(Y\to\infty\), equivalently as \(X\to\infty\),
\[
N_n(Y)=
\begin{cases}
\delta_n X+O_n(X^{1/2}), & n\ \text{odd},\\[4pt]
2\delta_n X+O_n(X^{1/2}), & n\ \text{even}.
\end{cases}
\]
\end{theorem}

\subsection{Related results}
There is an extensive literature on monogeneity, integral bases, and pure fields. For comprehensive treatments, we refer to the monographs of Evertse--Gy{\H{o}}ry~\cite{EG17} and Ga\'al~\cite{Gaal2019Book}. See also the recent surveys by Ga\'al~\cite{GaalSurvey2024} and Evertse--Gy{\H{o}}ry~\cite{EG24} for updated accounts of current developments.

More generally, for radical extensions \(L(\sqrt[n]{\alpha})/L\), Smith proved in \cite{SmithRadical2021} a relative criterion for \(\sqrt[n]{\alpha}\) to generate a power integral basis. Specializing to \(L=\Q\) and to the pure family \(X^n-m\) gives exactly the two local conditions in Theorem~\ref{thm:main}: square-freeness of \(m\) and \(\nu_p(m^p-m)=1\) for every \(p\mid n\). Earlier formulations for \(X^{p^r}-a\) and the associated Wieferich-type congruence \(a^p\equiv a\pmod {p^2}\) appear in \cite{Gassert2017,ElFadil2021Note}.

A companion paper \cite{ND-PeriodicShape} by the first-named author studies the full \(p\)-adic shape of explicit integral bases for pure fields. In that setting, if \(p^e\parallel n\), the local shape is determined by the class of the parameter modulo \(p^{e+1}\); equivalently, the global period is \(\prod_{p^e\parallel n}p^{e+1}=n\rad(n)\). The present paper uses only the coarser modulo \(p^2\) information needed to decide whether \(\Z[\sqrt[n]{m}]\) is maximal.

The density of square-free integers in arithmetic progressions and basic \(p\)-adic features of Frobenius, including Teichm\"uller lifts, enter our distributional results; see Brown's paper~\cite{Brown2021} and standard sources such as Cohen's book~\cite{CohenANTI} or Koblitz's book~\cite{Koblitzpadic}. For pure cubics, monogeneity and quantitative aspects were studied by Aygin--Nguyen~\cite{AyginNguyen2021,AyginNguyenCorrigendum2023}. In the monic coefficient space, Bhargava--Shankar--Wang~\cite{BSWI,BSWII} proved that a positive density \(\zeta(2)^{-1}\) of polynomials \(f\) satisfy \(\Z[x]/(f)=\mathcal O_{\Q[x]/(f)}\), and they determined the density of squarefree discriminants. Our work instead treats the rigid family \(f_m(X)=X^n-m\): we give a complete local-global criterion for \(\Z[\sqrt[n]{m}]\) to be maximal and compute the exact density \(\frac{6}{\pi^2}\prod_{p\mid n}\frac{p}{p+1}\). From a moduli-theoretic perspective, schemes parameterizing monogenic generators and notions of local monogenicity are developed in~\cite{ArpinBozleeSmith2023}; these results are not used here.

\subsection{Organization}
Section~\ref{sec:criterion} proves the criterion for \(\alpha\)-monogenic pure fields. Section~\ref{sec:density} develops the density calculation, the arithmetic-progression refinement, conditional-density independence of the local obstruction sets, and the discriminant-ordered field count.

\subsection{Acknowledgement}
The first author thanks the Morningside Center of Mathematics, Chinese Academy of Sciences, for its support and a stimulating research environment. We thank Professor Gy{\H{o}}ry for his interest in our work and for drawing our attention to his monograph with Evertse, as well as to the updated version of \cite{EG24}. We also thank Professor Ga\'al for several helpful communications. We appreciate the anonymous referee for the careful reading of our manuscript and for the detailed comments.

\section{A criterion for \texorpdfstring{$\alpha$}{alpha}-monogenic pure fields}\label{sec:criterion}

We recall some basic definitions and notations.
\begin{definition}
An algebraic number field \(K\) is called a \emph{pure number field of degree \(n\)} if \(K=\Q(\alpha)\), where \(\alpha\) is a root of an irreducible polynomial \(X^n-m\in\Z[X]\) with \(n\ge2\).
\end{definition}

\begin{definition}
Let \(K=\Q(\alpha)\), where \(\alpha\) is an algebraic integer. We say that \(K\) is \emph{\(\alpha\)-monogenic} if \(\Z[\alpha]=\OK\).
\end{definition}

For example, when \(K=\Q(\alpha)\) with \(\alpha\in\{\sqrt{-2},\sqrt[3]{3},\sqrt[5]{6}\}\), one checks that \(\OK=\Z[\alpha]\). The local computations below require only one standard tool: Dedekind's criterion for detecting when a prime divides the index of a monogenic order. We use the following form.

\begin{theorem}[Dedekind's index criterion]\label{thm:dedekind-index}
Let \(K=\Q(\alpha)\), where \(\alpha\) is integral over \(\Z\) with monic minimal polynomial \(f\in\Z[X]\). For a prime \(p\), write the factorization of \(\overline f\) in \(\F_p[X]\) as
\[
\overline f(X)=\overline\pi_1(X)^{e_1}\cdots \overline\pi_g(X)^{e_g},
\]
where the \(\overline\pi_j\) are pairwise distinct monic irreducible polynomials. Let \(\pi_j\in\Z[X]\) be monic lifts of \(\overline\pi_j\), and define \(F\in\Z[X]\) by
\[
f(X)=\pi_1(X)^{e_1}\cdots\pi_g(X)^{e_g}+pF(X).
\]
Then \(p\mid [\OK:\Z[\alpha]]\) if and only if \(\overline\pi_j\mid \overline F\) in \(\F_p[X]\) for some \(j\) with \(e_j\ge2\).
\end{theorem}

To restrict the possible prime divisors of the index, we also recall the discriminant-index identity.

\begin{theorem}\label{thm:index-disc}
Let \(K\) be a number field, and let \(\gamma_1,\dots,\gamma_d\) be a \(\Z\)-basis of a sublattice \(N\subset\OK\) of finite index. Then
\[
D(\gamma_1,\dots,\gamma_d)=[\OK:N]^2\disc(\OK/\Z).
\]
In particular, if \(K=\Q(\alpha)\), where \(\alpha\) is integral over \(\Z\) with minimal polynomial \(f\), then
\begin{equation}\label{eq:index-disc-special}
\disc(f)=[\OK:\Z[\alpha]]^2\disc(\OK/\Z).
\end{equation}
\end{theorem}

For the particular polynomial \(X^n-m\), the discriminant-index identity immediately localizes all possible index primes. This observation reduces maximality to checking only primes dividing \(mn\).

\begin{lemma}\label{lem:support}
Let \(K=\Q(\alpha)\), where \(\alpha\) has minimal polynomial \(f(X)=X^n-m\). Then every prime dividing \([\OK:\Z[\alpha]]\) divides \(mn\). Equivalently, \(K\) is \(\alpha\)-monogenic if and only if no prime \(p\mid mn\) divides \([\OK:\Z[\alpha]]\).
\end{lemma}

\begin{proof}
The polynomial discriminant is
\[
\disc(X^n-m)=(-1)^{(n-1)(n-2)/2}n^n m^{n-1}.
\]
The claim follows from \eqref{eq:index-disc-special}.
\end{proof}

In the case \(p\mid n\), the reduction modulo \(p\) naturally produces the expression \(m^{p^r}-m\), where \(p^r\) is the exact \(p\)-power contribution to \(n\). The following elementary lemma shows that, for the index criterion, this is equivalent to the simpler expression \(m^p-m\).

\begin{lemma}[Reduction from \(p^r\) to \(p\)]\label{lem:valuation-reduction}
Let \(p\) be a prime and let \(r\ge1\). If \(p\nmid m\), then
\[
\vp(m^{p^r}-m)=\vp(m^p-m).
\]
If \(p\mid m\) and \(\vp(m)=1\), then
\[
\vp(m^{p^r}-m)=\vp(m^p-m)=1.
\]
\end{lemma}

\begin{proof}
Assume first that \(p\nmid m\). Then
\[
\vp(m^{p^r}-m)=\vp(m^{p^r-1}-1).
\]
Write \(p^r-1=(p-1)u\), where \(u=1+p+\cdots+p^{r-1}\) is prime to \(p\). Put \(A=m^{p-1}\). If \(A=1\), then both sides below are \(+\infty\), and the asserted equality is immediate. Otherwise, since \(A\equiv1\pmod p\), the standard lifting-the-exponent formula gives
\[
\vp(A^u-1)=\vp(A-1)+\vp(u)=\vp(A-1).
\]
For \(p=2\), this is the elementary fact that \(\vp(x^u-1)=\vp(x-1)\) for odd \(x\) and odd \(u\), with the same harmless convention when \(x=1\). Hence
\[
\vp(m^{p^r}-m)=\vp(m^{p-1}-1)=\vp(m^p-m).
\]
If \(p\mid m\) and \(\vp(m)=1\), then
\[
m^{p^r}-m=m(m^{p^r-1}-1),\qquad m^p-m=m(m^{p-1}-1),
\]
and both factors in parentheses are congruent to \(-1\pmod p\). Thus both valuations are equal to \(1\).
\end{proof}

With these preparations, the proof of the criterion becomes a direct application of Dedekind's criterion at the finite set of possible index primes.
\begin{theorem}\label{thm:main}
Let \(K=\Q(\alpha)\), where \(\alpha\) is a root of the irreducible polynomial \(f(X)=X^n-m\in\Z[X]\). Then \(K\) is \(\alpha\)-monogenic, that is \(\Z[\alpha]=\OK\), if and only if
\begin{enumerate}[label=\textup{(\roman*)}]
\item \(m\) is square-free, and
\item \(\nu_p(m^p-m)=1\) for every prime \(p\mid n\).
\end{enumerate}
\end{theorem}

\begin{proof}
Let \(I=[\OK:\Z[\alpha]]\). By Lemma~\ref{lem:support}, only primes dividing \(mn\) can divide \(I\).

\smallskip
\noindent\emph{Primes dividing \(m\).}
Let \(p\mid m\). Reducing \(f\) modulo \(p\) gives \(\overline f(X)=X^n\). In Theorem~\ref{thm:dedekind-index}, take \(\pi(X)=X\). Then
\[
f(X)=X^n+pF(X),\qquad F(X)=-m/p\in\Z[X].
\]
Thus \(p\mid I\) if and only if \(X\mid \overline F\) in \(\F_p[X]\), which is equivalent to \(F\equiv0\pmod p\), or \(\vp(m)\ge2\). Therefore no prime divisor of \(m\) divides \(I\) if and only if \(m\) is square-free.

\smallskip
\noindent\emph{Primes dividing \(n\) but not \(m\).}
Let \(p\mid n\) and \(p\nmid m\). Write \(n=p^rs\) with \((p,s)=1\), and set
\[
g(X)=X^s-m.
\]
Since \(p\nmid sm\), the polynomial \(\overline g\in\F_p[X]\) is separable. Factor it as
\[
\overline g=\overline\pi_1\cdots\overline\pi_t
\]
with distinct monic irreducible factors, and choose monic lifts \(\pi_i\in\Z[X]\). Put \(H=\prod_i\pi_i\). Then \(H\equiv g\pmod p\), and
\[
\overline f(X)=X^{p^rs}-\overline m=(X^s-\overline m)^{p^r}=\overline g(X)^{p^r}=\prod_i\overline\pi_i(X)^{p^r}.
\]
Dedekind's criterion applies to
\[
\Phi(X)=\frac{f(X)-H(X)^{p^r}}{p}\in\Z[X].
\]
Because \(H\equiv g\pmod p\), we have \(H^{p^r}\equiv g^{p^r}\pmod {p^2}\). Indeed, if \(H=g+pB\), then every non-leading binomial term in \((g+pB)^{p^r}-g^{p^r}\) is divisible by \(p^2\). Hence
\begin{equation}\label{eq:Phi-replace-H-by-g}
\overline\Phi\equiv \overline{\frac{f-g^{p^r}}{p}}\pmod p.
\end{equation}
We next compute this polynomial modulo \(\overline g\). Since \(X^s\equiv m\pmod {g}\), the polynomial
\[
R(X)=f(X)-g(X)^{p^r}-(m^{p^r}-m)
\]
lies in the ideal \((g)\subset\Z[X]\). Moreover, \(R\) is coefficientwise divisible by \(p\): after expanding \((X^s-m)^{p^r}\), all intermediate binomial coefficients are divisible by \(p\), and the remaining constant term is also divisible by \(p\). Indeed, for \(p\) odd this constant term is \(0\), while for \(p=2\) it is
\(-2m^{2^r}\). Write \(R=gQ\). Reducing modulo \(p\) gives \(\overline g\,\overline Q=0\) in the domain \(\F_p[X]\), so \(\overline Q=0\), and hence \(Q=pQ_1\) with \(Q_1\in\Z[X]\). Therefore
\[
\frac{f-g^{p^r}}{p}-\frac{m^{p^r}-m}{p}=gQ_1,
\]
and so, in \(\F_p[X]/(\overline g)\),
\begin{equation}\label{eq:F-constant-mod-g}
\overline\Phi\equiv \overline{\frac{m^{p^r}-m}{p}}.
\end{equation}
The integer \((m^{p^r}-m)/p\) is well-defined because \(m^{p^r}\equiv m\pmod p\).

By \eqref{eq:F-constant-mod-g}, for every irreducible factor \(\overline\pi_i\) of \(\overline g\), the image of \(\overline\Phi\) in \(\F_p[X]/(\overline\pi_i)\) is the same constant \(\overline{(m^{p^r}-m)/p}\). Thus \(\overline\pi_i\mid\overline\Phi\) for some \(i\) if and only if
\[
\frac{m^{p^r}-m}{p}\equiv0\pmod p,
\]
that is, if and only if \(p^2\mid m^{p^r}-m\). Dedekind's criterion therefore gives
\[
p\mid I \iff \nu_p(m^{p^r}-m)\ge2.
\]
By Lemma~\ref{lem:valuation-reduction}, this is equivalent, under the present assumption \(p\nmid m\), to \(\nu_p(m^p-m)\ge2\).

Combining the two cases, \(I=1\) if and only if \(m\) is square-free and \(\nu_p(m^p-m)=1\) for every \(p\mid n\). Indeed, if \(p\mid m\), \(p\mid n\), and \(m\) is square-free at \(p\), Lemma~\ref{lem:valuation-reduction} gives \(\nu_p(m^p-m)=1\), so the condition at such a prime is automatic once square-freeness holds.
\end{proof}

The proof actually gives a little more than the vanishing criterion: it identifies exactly which primes can occur in the index.

\begin{corollary}\label{cor:prime-support}
With notation as in Theorem~\ref{thm:main}, the prime support of the index is
\[
\{p:p\mid[\OK:\Z[\alpha]]\}
=
\{p:\nu_p(m)\ge2\}\cup\{p\mid n:\nu_p(m^p-m)\ge2\}.
\]
\end{corollary}

\begin{proof}
The proof of Theorem~\ref{thm:main} shows that a prime divisor \(p\) of \(m\) divides the index exactly when \(\nu_p(m)\ge2\). It also shows that a prime \(p\mid n\) with \(p\nmid m\) divides the index exactly when \(\nu_p(m^p-m)\ge2\). If \(p\mid m\) and \(\nu_p(m)\ge2\), then
\[
\nu_p(m^p-m)=\nu_p\bigl(m(m^{p-1}-1)\bigr)=\nu_p(m),
\]
since \(m^{p-1}-1\equiv -1\pmod p\). Hence the displayed union gives exactly the prime divisors of the index.
\end{proof}

Specializing the criterion gives concrete congruence tests in small degrees.

\begin{corollary}\label{cor:small-degrees}
Let \(K=\Q(\alpha)\), where \(\alpha^n=m\) and \(X^n-m\) is irreducible over \(\Q\).
\begin{enumerate}[label=\textup{(\arabic*)}]
\item If \(n=2\), then \(K\) is \(\alpha\)-monogenic if and only if \(m\) is square-free and \(m\not\equiv1\pmod4\).
\item If \(n=3\), then \(K\) is \(\alpha\)-monogenic if and only if \(m\) is square-free and \(m\not\equiv\pm1\pmod9\).
\item If \(n=4\), then \(K\) is \(\alpha\)-monogenic if and only if \(m\) is square-free and \(m\not\equiv1\pmod4\).
\item If \(n=5\), then \(K\) is \(\alpha\)-monogenic if and only if \(m\) is square-free and \(m\not\equiv1,7,18,24\pmod{25}\).
\item If \(n=6\), then \(K\) is \(\alpha\)-monogenic if and only if \(m\) is square-free, \(m\not\equiv1\pmod4\), and \(m\not\equiv\pm1\pmod9\).
\end{enumerate}
\end{corollary}

\begin{proof}
Apply Theorem~\ref{thm:main}. For \(p=2\), the unique unit class modulo \(4\) satisfying \(x^2\equiv x\pmod4\) is \(1\), while square-free even \(m\) automatically has \(\nu_2(m^2-m)=1\). For \(p=3\), the excluded unit classes modulo \(9\) are \(\pm1\). For \(p=5\), the excluded unit classes are the Teichm\"uller lifts in \((\Z/25\Z)^\times\), namely \(1,7,18,24\). Combining the relevant primes gives the five assertions.
\end{proof}

\section{Natural density for \texorpdfstring{$\alpha$}{alpha}-monogenic pure fields}\label{sec:density}

Throughout this section, \(n\ge2\) is fixed.

\subsection*{Notation and preliminaries}
For \(q\ge1\) and \(a\in\Z\), let
\[
\A(q,a)=\{m\in\Z:m\equiv a\pmod q\}.
\]
For a set \(E\subset\Z\), we use the two-sided natural density
\[
d(E)=\lim_{X\to\infty}\frac{1}{2X}\#\{m\in E:0<|m|\le X\},
\]
provided the limit exists. If \(E^+\subset\Z_{>0}\), we write
\[
d_+(E^+)=\lim_{X\to\infty}\frac{1}{X}\#\{m\in E^+:1\le m\le X\}.
\]
All densities below are two-sided unless a one-sided density is explicitly indicated. Let \(\mu\) denote the M\"obius function; for \(m\ne0\), \(\mu^2(|m|)\) is the indicator of square-freeness. We record the elementary counting lemmas used below.

\begin{lemma}\label{lem:sf-progression}
Let \(q\ge1\) and let \(a\in\Z\) with \(\gcd(a,q)=1\). Then
\[
d\bigl(\{m\in\Z:m\equiv a\bmod q,\ m\text{ square-free}\}\bigr)
=\frac{1}{q\zeta(2)}\prod_{\ell\mid q}\frac{1}{1-\ell^{-2}}.
\]
\end{lemma}

\begin{proof}
Using \(\mu^2(|m|)=\sum_{d^2\mid m}\mu(d)\), we get
\[
\#\{0<|m|\le X:m\equiv a\bmod q,\ \mu^2(|m|)=1\}
=\sum_{d\le\sqrt X}\mu(d)N_d(X),
\]
where
\[
N_d(X)=\#\{0<|m|\le X:m\equiv a\bmod q,\ d^2\mid m\}.
\]
If \(\gcd(d,q)>1\), then \(N_d(X)=0\), because \(\gcd(a,q)=1\). If \(\gcd(d,q)=1\), the two congruence conditions combine, by the Chinese remainder theorem, into one residue class modulo \(qd^2\), and therefore
\[
N_d(X)=\frac{2X}{qd^2}+O(1).
\]
Extending the absolutely convergent main sum from \(d\le\sqrt X\) to all \(d\ge1\) changes the count by \(O_q(X^{1/2})\). It follows that
\[
\#\{0<|m|\le X:m\equiv a\bmod q,\ \mu^2(|m|)=1\}
=\frac{2X}{q}\sum_{\gcd(d,q)=1}\frac{\mu(d)}{d^2}+O_q(X^{1/2})
\qquad (X\to\infty).
\]
Since
\[
\sum_{\gcd(d,q)=1}\frac{\mu(d)}{d^2}
=\prod_{\ell\nmid q}(1-\ell^{-2})
=\frac{1}{\zeta(2)}\prod_{\ell\mid q}\frac{1}{1-\ell^{-2}},
\]
division by \(2X\) and passage to the limit prove the formula.
\end{proof}

When a new prime \(p\nmid q\) is introduced, we need to know how square-free integers distribute among the nonzero classes modulo \(p^2\) above a fixed reduced class modulo \(q\). The next lemma records the required uniformity.

\begin{lemma}[Adding a new \(p^2\)-coordinate]\label{lem:sf-new-p2}
Let \(q\ge1\), let \(a\in\Z\) with \((a,q)=1\), and let \(p\nmid q\) be prime. If \(b\in\Z/p^2\Z\) is nonzero, then the square-free integers satisfying
\[
m\equiv a\pmod q,\qquad m\equiv b\pmod {p^2}
\]
have two-sided density
\[
\frac{1}{p^2-1}\cdot
\frac{1}{q\zeta(2)}
\prod_{\ell\mid q}\frac{1}{1-\ell^{-2}}.
\]
\end{lemma}

\begin{proof}
If \(p\nmid b\), the two congruences combine to a class modulo \(qp^2\) which is coprime to \(qp^2\). Lemma~\ref{lem:sf-progression} gives the density
\[
\frac{1}{qp^2\zeta(2)}
\prod_{\ell\mid q}\frac{1}{1-\ell^{-2}}\frac{1}{1-p^{-2}}
=
\frac{1}{p^2-1}\cdot
\frac{1}{q\zeta(2)}
\prod_{\ell\mid q}\frac{1}{1-\ell^{-2}}.
\]
If \(p\mid b\) but \(p^2\nmid b\), write \(b=pb_1\) with \(b_1\not\equiv0\pmod p\). Then \(m=pu\), and the two congruences are equivalent to
\[
u\equiv p^{-1}a\pmod q,\qquad u\equiv b_1\pmod p.
\]
Moreover, \(m\) is square-free if and only if \(u\) is square-free; the condition \(p\nmid u\) is already forced by \(u\equiv b_1\pmod p\). Applying Lemma~\ref{lem:sf-progression} to \(u\) modulo \(qp\), and then replacing \(X\) by \(X/p\), gives the same density as above. 
\end{proof}

The residue classes excluded by the local criterion are the fixed points of the Frobenius map \(x\mapsto x^p\) modulo \(p^2\). Their structure is simple and will be used repeatedly.

\begin{lemma}[Fixed points of Frobenius modulo \(p^2\)]\label{lem:fixedpoints}
For a prime \(p\), the congruence
\[
x^p\equiv x\pmod {p^2}
\]
has exactly \(p\) solutions modulo \(p^2\). One solution is \(0\), and the other \(p-1\) are the Teichm\"uller lifts in \((\Z/p^2\Z)^\times\). Moreover, each solution modulo \(p\) lifts uniquely to a solution modulo \(p^k\) for every \(k\ge1\).
\end{lemma}

\begin{proof}
Let \(h(x)=x^p-x\). In \(\F_p\), every element is a root of \(h\), and at each such root one has
\[
h'(x)=px^{p-1}-1\equiv -1\pmod p,
\]
which is invertible. Hensel's lemma gives a unique lift modulo \(p^k\) for every \(k\ge1\). In particular, there are exactly \(p\) solutions modulo \(p^2\), one above each class modulo \(p\). The nonzero lifted classes are precisely the Teichm\"uller representatives of \(\F_p^\times\).
\end{proof}

\subsection{Local obstruction classes}

For a prime \(p\), define
\[
E_p=\{a\in\Z/p^2\Z:a^p\equiv a\pmod {p^2}\},
\qquad
\E_p=\{m\in\Z:m\bmod p^2\in E_p\}.
\]
By Lemma~\ref{lem:fixedpoints}, \(E_p\) consists of the class \(0\bmod p^2\) and the \(p-1\) Teichm\"uller unit classes. For square-free \(m\), the class \(0\bmod p^2\) cannot occur. By Theorem~\ref{thm:main}, if \(p\mid n\), then \(\E_p\) is exactly the set on which the \(p\)-local condition \(\nu_p(m^p-m)=1\) fails.

Let \(\Ssq\) denote the set of nonzero square-free integers. We now translate the congruence \(m^p\equiv m\pmod {p^2}\) into a relative density inside \(\Ssq\).

\begin{proposition}[One-prime local density]\label{prop:local-density}
For every prime \(p\),
\[
\frac{d(\Ssq\cap\E_p)}{d(\Ssq)}=\frac{1}{p+1},
\qquad
\frac{d(\Ssq\setminus\E_p)}{d(\Ssq)}=\frac{p}{p+1}.
\]
\end{proposition}

\begin{proof}
The square-free integers do not meet the class \(0\bmod p^2\). Among the unit classes modulo \(p^2\), exactly \(p-1\) classes belong to \(E_p\). By Lemma~\ref{lem:sf-progression}, each such unit class contributes two-sided density
\[
\frac{1}{p^2\zeta(2)}\frac{1}{1-p^{-2}}
=\frac{1}{\zeta(2)}\frac{1}{p^2-1}.
\]
Consequently
\[
d(\Ssq\cap\E_p)=\frac{1}{\zeta(2)}\frac{p-1}{p^2-1}
=\frac{1}{\zeta(2)}\frac{1}{p+1}.
\]
Since \(d(\Ssq)=1/\zeta(2)\), the two asserted relative densities follow.
\end{proof}

\subsection{Global density}

For a finite set \(P\) of primes, put
\[
\mathcal B_P=\Ssq\cap\bigcap_{p\in P}\E_p.
\]
The set \(\bigcap_{p\in P}\E_p\), after imposing square-freeness, is a disjoint union of \(\prod_{p\in P}(p-1)\) unit residue classes modulo \(\prod_{p\in P}p^2\). By the Chinese remainder theorem and Lemma~\ref{lem:sf-progression},
\begin{equation}\label{eq:BP-density}
d(\mathcal B_P)=\frac{1}{\zeta(2)}\prod_{p\in P}\frac{p-1}{p^2-1}
=\frac{1}{\zeta(2)}\prod_{p\in P}\frac{1}{p+1}.
\end{equation}
Now set
\[
\mathcal G_n=\Ssq\cap\bigcap_{p\mid n}(\Z\setminus\E_p).
\]
Inclusion-exclusion gives
\begin{equation}\label{eq:Gn-density}
d(\mathcal G_n)=\frac{1}{\zeta(2)}\prod_{p\mid n}\left(1-\frac{1}{p+1}\right)
=\frac{1}{\zeta(2)}\prod_{p\mid n}\frac{p}{p+1}.
\end{equation}

The density computation is made among all nonzero integer parameters, while the field \(K_m\) has degree \(n\) only when \(X^n-m\) is irreducible. The next lemma shows that this restriction does not affect the main term.

\begin{lemma}\label{lem:irr-density}
The set
\[
\mathcal I_n=\{m\in\Z\setminus\{0\}:X^n-m\text{ is irreducible over }\Q\}
\]
has two-sided natural density \(1\). More precisely,
\[
\#\{0<|m|\le X:m\notin\mathcal I_n\}=O_n(X^{1/2})\qquad (X\to\infty).
\]
\end{lemma}

\begin{proof}
If \(m\) has a prime divisor \(p\) with \(\nu_p(m)=1\), then \(X^n-m\) is Eisenstein at \(p\), hence irreducible over \(\Q\). Thus the complement of \(\mathcal I_n\) is contained in \(\{\pm1\}\) together with the squarefull integers, namely the nonzero integers \(m\) such that every prime exponent in \(|m|\) is at least \(2\). The number of squarefull integers of absolute value at most \(X\) is \(O(X^{1/2})\): each positive squarefull integer can be written as \(a^2b^3\) with \(b\) square-free, and
\[
\sum_{b\le X^{1/3}}\left(\frac{X}{b^3}\right)^{1/2}=O(X^{1/2}).
\]
The same bound holds after allowing both signs.
\end{proof}

Combining the local density calculation with inclusion-exclusion over the finitely many primes dividing \(n\), we obtain the global density.
\begin{theorem}\label{thm:density}
Let
\[
\M_n=\{m\in\Z\setminus\{0\}:X^n-m\text{ is irreducible over }\Q\text{ and }\Z[\alpha_m]=\mathcal O_{K_m}\}.
\]
Then
\[
d(\M_n)=\delta_n:=\frac{6}{\pi^2}\prod_{p\mid n}\frac{p}{p+1}.
\]
In particular, \(\delta_n\) depends only on \(\rad(n)\).
\end{theorem}

\begin{proof}
By Theorem~\ref{thm:main}, for \(m\in\mathcal I_n\) the equality \(\Z[\alpha_m]=\mathcal O_{K_m}\) holds if and only if \(m\in\mathcal G_n\). Hence
\[
\M_n=\mathcal I_n\cap\mathcal G_n.
\]
By Lemma~\ref{lem:irr-density}, the complement of \(\mathcal I_n\) has density zero and contributes \(O_n(X^{1/2})\) points with \(0<|m|\le X\). Therefore \(\M_n\) has the same density as \(\mathcal G_n\). Formula \eqref{eq:Gn-density} gives
\[
d(\M_n)=\frac{1}{\zeta(2)}\prod_{p\mid n}\frac{p}{p+1}
=\frac{6}{\pi^2}\prod_{p\mid n}\frac{p}{p+1}.
\]
Finally, the product depends only on the primes dividing \(n\), hence only on \(\rad(n)\).
\end{proof}

\begin{example}[Degree-specific densities]\label{ex:degree-specific-density}
For \(n=4\), Theorem~\ref{thm:main} gives
\[
\OK=\Z[\alpha]\iff m\text{ is square-free and }m\not\equiv1\pmod4.
\]
For \(n=6\), it gives
\[
\OK=\Z[\alpha]
\iff
m\text{ is square-free},\quad m\not\equiv1\pmod4,
\quad m\not\equiv\pm1\pmod9.
\]
The corresponding densities are
\[
\delta_4=\frac{4}{\pi^2},\qquad \delta_6=\frac{3}{\pi^2}.
\]
\end{example}

\begin{proof}
For \(n=4\), the only prime dividing \(n\) is \(2\). If \(m\) is square-free and even, then \(\nu_2(m^2-m)=\nu_2(m)=1\). If \(m\) is odd, then \(\nu_2(m^2-m)=\nu_2(m-1)\), which equals \(1\) exactly when \(m\equiv3\pmod4\). Thus the excluded square-free class is \(m\equiv1\pmod4\).

For \(n=6\), the prime \(2\) gives the same condition. For \(p=3\), since \(m^3-m=m(m-1)(m+1)\), a square-free \(m\) divisible by \(3\) automatically satisfies \(\nu_3(m^3-m)=1\), while a unit \(m\) modulo \(3\) satisfies \(\nu_3(m^3-m)\ge2\) exactly for \(m\equiv\pm1\pmod9\). The density values follow from Theorem~\ref{thm:density}.
\end{proof}

\begin{remark}[Pure cubics and field monogeneity]\label{rem:cubic}
For square-free \(m\), Theorem~\ref{thm:main} gives
\[
\mathcal O_{\Q(\sqrt[3]{m})}=\Z[\sqrt[3]{m}]
\iff
m\not\equiv\pm1\pmod9.
\]
If \(m\equiv\pm1\pmod9\), then \(\Z[\sqrt[3]{m}]\) is not the maximal order. The field may still be monogenic only if some algebraic integer other than the distinguished root generates the full ring of integers. This is the additional global question studied for pure cubics in \cite{AyginNguyen2021}; the comparison quoted in Corollary~\ref{cor:intro-comparison} says that the exceptional square-free parameters in these two congruence classes have density zero. 
\end{remark}

\subsection{Consequences and applications of the density formula}

Write \(\delta_n=\frac{6}{\pi^2}\prod_{p\mid n}\frac{p}{p+1}\). The same proof gives an error term for the associated counting function.

\begin{corollary}\label{cor:interval-count}
Fix \(n\ge2\). Let
\[
\mathcal M_n(X)=\#\{m\in\M_n:0<|m|\le X\}.
\]
Then, as \(X\to\infty\),
\[
\mathcal M_n(X)=2\delta_nX+O_n(X^{1/2}).
\]
If one restricts to \(1\le m\le X\), the corresponding count is \(\delta_nX+O_n(X^{1/2})\) as \(X\to\infty\).
\end{corollary}

\begin{proof}
The error term follows from the inclusion-exclusion proof of
\eqref{eq:Gn-density}: for each \(P\subseteq\{p:p\mid n\}\), the set
\(\Ssq\cap\bigcap_{p\in P}\E_p\) is, after square-freeness removes the
zero classes, a finite union of unit residue classes modulo
\(\prod_{p\in P}p^2\). Lemma~\ref{lem:sf-progression} gives an
\(O_P(X^{1/2})\) error for each such class. Summing over the finitely many
\(P\) and using Lemma~\ref{lem:irr-density} gives the asserted
\(O_n(X^{1/2})\) error. 
\end{proof}

\subsection{Arithmetic progressions}

For a prime \(p\), recall that
\[
E_p=\{a\in\Z/p^2\Z:a^p\equiv a\pmod {p^2}\}.
\]
The argument is local enough to keep track of an externally imposed arithmetic progression. The following formula records how the \(p\)-local factor changes according to the \(p\)-adic depth of the modulus.

\begin{corollary}[Density in arithmetic progressions]\label{cor:AP-density}
Let \(q\ge1\) and \(a\in\Z\) with \((a,q)=1\). Define
\[
\M_n(q,a)=\{m\in\Z:m\equiv a\pmod q,
\ X^n-m\text{ is irreducible over }\Q,
\ \Z[\alpha_m]=\mathcal O_{K_m}\}.
\]
Then \(\M_n(q,a)\) has two-sided natural density
\begin{equation}\label{eq:AP-density}
d(\M_n(q,a))=
\frac{1}{q\zeta(2)}
\prod_{\ell\mid q}\frac{1}{1-\ell^{-2}}
\prod_{p\mid n}\lambda_p(q,a),
\end{equation}
where
\[
\lambda_p(q,a)=
\begin{cases}
\displaystyle \frac{p}{p+1}, & \nu_p(q)=0,\\[6pt]
\displaystyle \frac{p-1}{p}, & \nu_p(q)=1,\\[6pt]
1, & \nu_p(q)\ge2\text{ and }a\bmod p^2\notin E_p,\\[4pt]
0, & \nu_p(q)\ge2\text{ and }a\bmod p^2\in E_p.
\end{cases}
\]
The same formula is the one-sided density of positive \(m\) in the same progression. Equivalently, the relative density inside the progression \(m\equiv a\pmod q\) is
\[
\frac{1}{\zeta(2)}
\prod_{\ell\mid q}\frac{1}{1-\ell^{-2}}
\prod_{p\mid n}\lambda_p(q,a).
\]
\end{corollary}

\begin{proof}
By Lemma~\ref{lem:sf-progression}, square-free integers in the progression \(m\equiv a\pmod q\) have density
\[
\frac{1}{q\zeta(2)}\prod_{\ell\mid q}\frac{1}{1-\ell^{-2}}.
\]
It remains to impose, for each \(p\mid n\), the condition \(m\notin\E_p\).

If \(\nu_p(q)=0\), no residue information modulo \(p\) is prescribed. Refining the fixed class \(a\bmod q\) by classes modulo \(p^2\), Lemma~\ref{lem:sf-new-p2} shows that every nonzero class modulo \(p^2\) has the same square-free density. The class \(0\bmod p^2\) contributes no square-free integers, and among the remaining \(p^2-1\) classes exactly the \(p-1\) Teichm\"uller unit classes lie in \(E_p\). Thus the relative factor for \(m\notin\E_p\) is
\[
1-\frac{p-1}{p^2-1}=\frac{p}{p+1}.
\]
If \(\nu_p(q)=1\), then \(a\not\equiv0\pmod p\), and among the \(p\) unit classes modulo \(p^2\) above \(a\bmod p\), exactly one lies in \(E_p\). These \(p\) unit classes have equal square-free density by Lemma~\ref{lem:sf-progression}, so the relative factor is \((p-1)/p\). If \(\nu_p(q)\ge2\), the class modulo \(p^2\) is fixed, giving factor \(1\) or \(0\) according as that class does not lie, or lies, in \(E_p\). Combining the finitely many local conditions by the Chinese remainder theorem gives \eqref{eq:AP-density}. The irreducibility condition changes the count only by \(O_n(X^{1/2})\), by Lemma~\ref{lem:irr-density}.
\end{proof}

\begin{example}[Quartic case in a fixed progression]\label{ex:quartic-AP}
Let \(n=4\), \(q=4\), and \(a\equiv3\pmod4\). Then \(\nu_2(q)=2\) and \(a\bmod4\notin E_2=\{0,1\}\), so \(\lambda_2(q,a)=1\). Therefore
\[
d\bigl(\{m\in\Z:m\equiv3\pmod4,\ \Z[\sqrt[4]{m}]=\mathcal O_{\Q(\sqrt[4]{m})}\}\bigr)
=\frac{1}{4\zeta(2)}\frac{1}{1-2^{-2}}
=\frac{2}{\pi^2}.
\]
On this progression, the condition at \(p=2\) is automatic; only square-freeness remains.
\end{example}

At modulus \(p^2\), each local condition is either forced or ruled out by the chosen residue class. This gives progressions on which all local conditions are automatic.

\begin{corollary}[Progressions on which the local conditions are automatic]\label{cor:full-AP}
Let
\[
M=\prod_{p\mid n}p^2.
\]
There exist residue classes \(a\bmod M\) with \(\gcd(a,M)=1\) such that every square-free integer \(m\equiv a\pmod M\) satisfies \(\Z[\alpha_m]=\mathcal O_{K_m}\), provided \(X^n-m\) is irreducible. For each such class \(a\), the set of such \(m\) has density
\[
\frac{1}{M\zeta(2)}\prod_{\ell\mid M}\frac{1}{1-\ell^{-2}}.
\]
\end{corollary}

\begin{proof}
For each \(p\mid n\), choose a unit class \(a_p\in(\Z/p^2\Z)^\times\setminus E_p\). Such a class exists because \(E_p\) contains only \(p-1\) unit classes, while \((\Z/p^2\Z)^\times\) has \(p^2-p\) classes. By the Chinese remainder theorem, choose \(a\bmod M\) satisfying \(a\equiv a_p\pmod {p^2}\) for every \(p\mid n\). Then \(\gcd(a,M)=1\), and every square-free \(m\equiv a\pmod M\) satisfies \(m\notin\E_p\) for every \(p\mid n\). Theorem~\ref{thm:main} gives \(\Z[\alpha_m]=\mathcal O_{K_m}\) whenever \(X^n-m\) is irreducible. The density is the square-free density in the progression, given by Lemma~\ref{lem:sf-progression}.
\end{proof}

\subsection{Discriminant ordering and local obstruction patterns}

On the set of parameters \(m\) satisfying \(\Z[\alpha_m]=\mathcal O_{K_m}\), the index-discriminant identity turns the parameter height \(|m|\) into the field-discriminant height. We now make this conversion precise.

\begin{theorem}\label{thm:fields-by-disc}
Let \(N_n(Y)\) denote the number of isomorphism classes of fields \(K=\Q(\alpha)\) for which \(\alpha^n=m\) for some \(m\in\Z\setminus\{0\}\), the polynomial \(X^n-m\) is irreducible over \(\Q\), \(\Z[\alpha]=\mathcal O_K\), and \(|\disc(K)|\le Y\). Put
\[
X=\left(\frac{Y}{n^n}\right)^{1/(n-1)}.
\]
Then, as \(Y\to\infty\), equivalently as \(X\to\infty\),
\[
N_n(Y)=
\begin{cases}
\delta_nX+O_n(X^{1/2}), & n\ \text{odd},\\[4pt]
2\delta_nX+O_n(X^{1/2}), & n\ \text{even}.
\end{cases}
\]
\end{theorem}

\begin{proof}
If \(\Z[\alpha]=\mathcal O_K\), then \eqref{eq:index-disc-special} gives
\[
\disc(K)=\disc(X^n-m)=(-1)^{(n-1)(n-2)/2}n^n m^{n-1}.
\]
Thus \(|\disc(K)|\le Y\) is equivalent to \(|m|\le X\).

It remains only to pass from parameters to fields. For parameters satisfying \(\Z[\alpha]=\mathcal O_K\), the discriminant determines \(|m|\). If \(n\) is odd, the parameters \(m\) and \(-m\) define the same field, because replacing \(\alpha\) by \(-\alpha\) changes \(m\) to \(-m\). Thus it suffices to count positive \(m\), giving \(\delta_nX+O_n(X^{1/2})\) as \(X\to\infty\). If \(n\) is even, the discriminants attached to \(m\) and \(-m\) have opposite signs, so the corresponding fields are not isomorphic. Thus all signed parameters are counted, giving \(2\delta_nX+O_n(X^{1/2})\) as \(X\to\infty\) by Corollary~\ref{cor:interval-count}.
\end{proof}

\begin{example}\label{ex:monotonicity}
Since
\[
\delta_n=\frac{6}{\pi^2}\prod_{p\mid n}\frac{p}{p+1},
\]
adjoining a new prime divisor to \(n\) multiplies \(\delta_n\) by \(p/(p+1)<1\). Hence if \(n\mid n'\), then \(\delta_{n'}\le\delta_n\). For instance,
\[
\delta_2=\delta_4=\frac{4}{\pi^2}\approx0.4053,
\quad
\delta_3=\frac{9}{2\pi^2}\approx0.4559,
\quad
\delta_5=\frac{5}{\pi^2}\approx0.5066,
\quad
\delta_6=\frac{3}{\pi^2}\approx0.3040.
\]
\end{example}

The same finite-product computation also describes the full pattern of local obstructions, not only the pattern in which none occur.

\begin{proposition}[Local obstruction patterns]\label{prop:obstruction-patterns}
Let \(P_n=\{p:p\mid n\}\). For each subset \(T\subseteq P_n\), define
\[
\mathcal F_T=\bigl\{m\in\Ssq:\{p\in P_n:m\in\E_p\}=T\bigr\}.
\]
Then
\[
d(\mathcal F_T)=\frac{1}{\zeta(2)}
\prod_{p\in T}\frac{1}{p+1}
\prod_{p\in P_n\setminus T}\frac{p}{p+1}.
\]
Equivalently, the sets \(\E_p\), \(p\in P_n\), are independent relative to the ambient set \(\Ssq\) in the sense that, for every \(T\subseteq P_n\),
\[
\frac{d\bigl(\Ssq\cap\bigcap_{p\in T}\E_p\bigr)}{d(\Ssq)}
=\prod_{p\in T}\frac{d(\Ssq\cap\E_p)}{d(\Ssq)}.
\]
The same formulas remain valid after imposing the additional condition that \(X^n-m\) be irreducible over \(\Q\).
\end{proposition}

\begin{proof}
The displayed independence identity is exactly \eqref{eq:BP-density} divided by \(d(\Ssq)=1/\zeta(2)\), together with Proposition~\ref{prop:local-density}. The formula for the exact pattern \(\mathcal F_T\) follows by inclusion-exclusion, or equivalently by multiplying the relative factors \(1/(p+1)\) for \(p\in T\) and \(p/(p+1)\) for \(p\notin T\), and then multiplying by \(d(\Ssq)=1/\zeta(2)\). Lemma~\ref{lem:irr-density} shows that the reducible parameters have density zero and contribute only \(O_n(X^{1/2})\) points of height at most \(X\).
\end{proof}

Finally, we record a local conditional form of the same calculation. It explains how much of the \(p^2\)-obstruction remains after one has fixed only the residue class modulo \(p\).

\begin{proposition}[Conditioning on a residue class modulo \(p\)]\label{prop:conditional-local}
Let \(p\mid n\), and let \(a\in\Z/p\Z\). Put
\[
\Ssq_a=\{m\in\Ssq:m\equiv a\pmod p\}.
\]
Then the relative density inside \(\Ssq_a\) of the set satisfying the \(p\)-local condition \(m\notin\E_p\) is
\[
\frac{d(\Ssq_a\setminus\E_p)}{d(\Ssq_a)}=
\begin{cases}
1, & a\equiv0\pmod p,\\[4pt]
\displaystyle\frac{p-1}{p}, & a\not\equiv0\pmod p.
\end{cases}
\]
Thus imposing \(m\equiv0\pmod p\) makes the \(p\)-local condition automatic for square-free \(m\), while imposing any nonzero class modulo \(p\) replaces the unrestricted relative factor \(p/(p+1)\) by \((p-1)/p\).
\end{proposition}

\begin{proof}
If \(a\equiv0\pmod p\) and \(m\in\Ssq_a\), then \(p\mid m\) but \(p^2\nmid m\), so \(m\notin\E_p\); equivalently, \(\nu_p(m^p-m)=1\). If \(a\not\equiv0\pmod p\), then the \(p\) residue classes modulo \(p^2\) lying above \(a\) are all unit classes, and Lemma~\ref{lem:fixedpoints} says that exactly one of them lies in \(E_p\). By Lemma~\ref{lem:sf-progression}, these \(p\) unit classes have equal square-free density. Hence the relative density of classes outside \(E_p\) is \((p-1)/p\).
\end{proof}

\begin{remark}[From \(p\)-level to \(p^2\)-level conditioning]\label{rem:depth}
Proposition~\ref{prop:conditional-local} conditions only on a class modulo \(p\). If one instead fixes a class \(b\bmod p^2\), then, among square-free integers in that class, the \(p\)-local condition \(m\notin\E_p\) is either automatic or impossible: it is automatic when \(b\notin E_p\), and impossible when \(b\in E_p\). This is the mechanism behind Corollary~\ref{cor:full-AP}: choosing a modulus divisible by \(p^2\) for every \(p\mid n\) and choosing, at each such prime, a class outside \(E_p\), forces the corresponding local condition throughout the progression.
\end{remark}

\bibliographystyle{alpha}
\bibliography{References}

\end{document}